\documentclass{amsart} 

\usepackage{tikz}
\usepackage{tikz-cd}
\usepackage{amsmath}
\usepackage{amsfonts}
\usepackage{amssymb}
\usepackage{amsthm}
\usepackage{mathrsfs}
\usepackage{mathtools}
\usepackage{color}
\usepackage{graphicx}
\usepackage{comment}
\usetikzlibrary{automata,topaths}
\usepackage{smallgrf}
\usepackage{graphicx}  
\usepackage[all]{xy}
\usepackage{amscd}
\usepackage{hyperref}
\hypersetup{%
	citebordercolor={0 0 1},
	linkbordercolor={0 0.25 0.25},
	urlbordercolor={1 1 1}
}

\title{Phantom maps and rational homotopy}

\author{Hiroshi Kihara}
\address{Center for Mathematical Sciences, University of Aizu, 
	Tsuruga, Ikki-machi, Aizu-Wakamatsu City, Fukushima, 965-8580, Japan}
\email{(kihara@u-aizu.ac.jp)}


\newtheoremstyle{mytheoremstyle} 
{1.25\topsep}                
{2.0\topsep}                 
{\it}
{}                           
{}                  
{}                           
{1em}                        
{\thmname{#1}\ \thmnumber{#2}\  \thmnote{#3}\hspace{0pt}} 

\newtheorem{thm}{Theorem}[section]    
\newtheorem{lem}[thm]{Lemma}          
\newtheorem{prop}[thm]{Proposition}

\theoremstyle{definition}

\newtheorem{cor}[thm]{Corollary}

\theoremstyle{mytheoremstyle}

\newtheorem{exa}[thm]{Example}

\newtheorem{rem}[thm]{Remark}

\def\fbb{{\mathbb{F}}}
\def\gbb{{\mathbb{G}}}

\def\qbb{{\mathbb{Q}}}

\def\zbb{{\mathbb{Z}}}

\def\acal{{\mathcal{A}}}
\def\bcal{{\mathcal{B}}}
\def\ccal{{\mathcal{C}}}

\def\ncal{{\mathcal{N}}}

\def\qcal{{\mathcal{Q}}}

\def\wcal{{\mathcal{W}}}

\def\hom{{\mathrm{Hom}}}

\def\ker{{\mathrm{Ker}}}

\numberwithin{equation}{section}

\newcommand{\rmsph}{{\mathrm{SPh}}}
\newcommand{\rmph}{{\mathrm{Ph}}} 

\definecolor{my_color}{rgb}{0, 0, 0}

\newcommand{\rmitem}[1]{\item[{\rm{(#1)}}] }

\begin{document}

\begin{abstract}
	We generalize theorems of McGibbon-Roitberg, Iriye, and Meier on the relations between phantom maps and rational homotopy, and apply them to provide new calculational examples of the homotopy sets $\rmph(X,Y)$ of phantom maps and the subsets $\rmsph(X,Y)$ of special phantom maps.
\end{abstract}

\maketitle

\section{Introduction} \label{section1}

The concept of a phantom map is a key to understanding maps with
infinite dimensional sources, and has been an important topic in homotopy theory since its discovery (\cite{McGibbon95,Roitberg94}). In this paper, we investigate the relationship between phantom maps and rational homotopy, which is a central concern in the study of phantom maps (see \cite{McGibbon95, Roitberg94, MR94, MR94II, I02, Iriye04, Meier78}) \par
In this section, we make a brief review of the basic notions and results on phantom maps, and present the main results of this paper.\par
Given two pointed $CW$-complexes $X$ and $Y$, a map $f: X \longrightarrow Y$ is called a {\em{phantom map}} if for any finite complex $K$ and any map $h: K \longrightarrow X$, the composite $fh$ is null homotopic.
Let $\rmph(X, Y)$ denote the subset of $[X, Y]$ consisting of homotopy classes of phantom maps.\par
We are also interested in the subset $\rmsph(X,Y)$ of $\rmph(X,Y)$ consisting of homotopy classes of {\em{special phantom maps}}, defined by the exact sequence of pointed sets
\begin{equation}\label{1.1}
\nonumber 0 \longrightarrow \rmsph(X,Y) \longrightarrow \rmph(X,Y) \overset{e_{Y\sharp}}{\longrightarrow} \rmph(X,\check{Y}), \eqno(1. 1)
\end{equation}
where $e_Y: Y \longrightarrow \check{Y} = \underset{p}{\prod}\ Y_{(p)}$ is a natural map called the \textit{local expansion} (cf. \cite[p.~150]{Roitberg94}). The target $Y$ is usually assumed to be nilpotent of finite type.\par
The following theorem is due to McGibbon-Roitberg and Iriye; the statements on phantom maps and special phantom maps are Theorem 2 in \cite{MR94} and Theorem 1.4 in \cite{I02} respectively. See Remark \ref{rem1.1} for the basic notions used in Theorem A.\\
\par\noindent {\bf{Theorem A}}\label{A} {\it{%
		{\rm (1)} Let $f:X'\longrightarrow X$ be a map between finite type domains and $Y$ a finite type target. If $f$ induces a monomorphism on the rational homology, then the maps
\begin{eqnarray*}
	f^\sharp:\rmph(X,Y) \longrightarrow \rmph(X',Y),\\
	f^\sharp:\rmsph(X,Y) \longrightarrow \rmsph(X',Y)
\end{eqnarray*}
are surjective.\\
{\rm (2)} Let $X$ be a finite type domain and $g:Y\longrightarrow Y'$ a map between finite type targets. If $g$ induces an epimorphism on the rational homotopy groups of degree $\geq 2$, then the maps
\begin{eqnarray*}
	g_\sharp:\rmph(X,Y) \longrightarrow \rmph(X,Y'),\\
	g_\sharp:\rmsph(X,Y) \longrightarrow \rmsph(X,Y')
\end{eqnarray*}
are surjective.}}\vspace{2mm}
\begin{rem}\label{rem1.1}{\rm 
	McGibbon-Roitberg \cite{MR94} and Iriye \cite{I02} dealt with skeletally phantom maps under finite type conditions on domains and targets; recall that a connected $CW$-complex is called a {\sl finite type domain} (resp. {\sl finite type target}) if its integral homology groups (resp. homotopy groups) are finitely generated in each degree. (See \cite[p. 30]{no10} and \cite[Section 1]{McGibbon95} for the difference between the two notions of phantom maps.) However, the most important class of finite type domains is that of $CW$-complexes of finite type (i.e. $CW$-complexes with finite skeleta), for which the two notions of phantom maps coincide.
}\end{rem}\vspace{2mm}
The following theorem is due to Meier (\cite[Theorem 5]{Meier78}). Recall that a space is called an {\sl $H_0$-space} if its rationalization is homotopy equivalent to a product of Eilenberg-MacLane spaces. Let $\hat{\zbb}$ denote the product $\Pi_p \hat{\zbb}_p$ of the $p$-completions of $\zbb$, in which $\zbb$ is diagonally contained.\\
\par\noindent {\bf{Theorem B}}\label{B} {\it{%
Let $Y$ be a nilpotent $CW$-complex of finite type which is an $H_0$-space. Suppose that there exists a finite product $Y'$ of copies of $BU$ and $U$, and a map $g:Y\longrightarrow Y'$ inducing a monomorphism on the rational homotopy groups. Then,
\[
	\rmph(K(\zbb,n),Y)\cong \underset{i>0}{\prod} H^i(K(\zbb,n); \pi_{i+1}(Y)\otimes \hat{\zbb}/\zbb) 
\]
holds for $n\ge 3$.}}\vspace{2mm}\par
We give generalizations of Theorems A and B as two main theorems, and apply them to calculations of $\rmph(X,Y)$ and $\rmsph(X,Y)$. The proofs are based on the approach introduced in a previous paper \cite{phantom}, which largely extends the rationalization-completion approach developed by Meier and Zabrodsky (\cite[Section 5]{McGibbon95}). Our approach is so general and enables us to give a generalization of Theorem A along with its simple proof. We also provide a generalization of Theorem B along with its dual version; the important idea of Theorem B has been largely overlooked in the literature.\par
Let $\ccal\wcal$ denote the category of pointed connected $CW$-complexes and homotopy classes of maps and let $\ncal$ denote the full subcategory of $\ccal\wcal$ consisting of nilpotent $CW$-complexes of finite type.\par
The first main theorem is the following.
\begin{thm}\label{thm1.1}
	Let $(f^{op},g):(X,Y)\longrightarrow (X',Y')$ be a morphism of $\ccal\wcal^{op}\times \ncal$. If $H^i(f;\pi_{i+1}(g)\otimes \mathbb{Q})$ is an epimorphism for any $i>0$, then the maps
	\begin{eqnarray*}
		(f^\sharp,g_\sharp):\rmph(X,Y)\longrightarrow \rmph(X',Y'),\\
		(f^\sharp,g_\sharp):\rmsph(X,Y) \longrightarrow \rmsph(X',Y')
	\end{eqnarray*}
	are surjective.
\end{thm}
The following corollary, and hence Theorem \ref{thm1.1} can be regarded as a generalization of Theorem A. 
\begin{cor}\label{cor1.2}
	\begin{itemize}
	\item[$(1)$] Let $f:X'\longrightarrow X$ be a map of $\ccal\wcal$, and $Y$ an object of $\ncal$. If $H_i(f;\mathbb{Q})$ is a monomorphism for $i>0$ with $\pi_{i+1}(Y)\otimes \mathbb{Q} \neq 0$, then the maps
	\begin{eqnarray*}
		f^\sharp:\rmph(X,Y) \longrightarrow \rmph(X',Y),\\
		f^\sharp:\rmsph(X,Y) \longrightarrow \rmsph(X',Y)
	\end{eqnarray*}
	are surjective.
	\item[$(2)$] Let $X$ be an object of $\ccal\wcal$ and $g:Y\longrightarrow Y'$ a map of $\ncal$. If $\pi_{i+1}(g)\otimes \mathbb{Q}$ is an epimorphism for $i>0$ with $H_i(X;\mathbb{Q})\neq 0$, then the maps
	\begin{eqnarray*}
		g_\sharp:\rmph(X,Y) \longrightarrow \rmph(X,Y'),\\
		g_\sharp:\rmsph(X,Y) \longrightarrow \rmsph(X,Y')
	\end{eqnarray*}
	are surjective.
	\end{itemize}
\end{cor}\vspace{2mm}
Corollary \ref{cor1.2} is used to obtain new vanishing results for $\rmph(X,Y)$ (see Proposition \ref{prop1.4}, Corollary \ref{cor1.5}, and Example \ref{exa2.5}).\par
Next, to state the second main theorem, we recall the basics of $\rmph(X,Y)$ and $\rmsph(X,Y)$ from \cite{phantom}. The set $\rmph(X,Y)$ (resp. $\rmsph(X,Y)$) can be described as the orbit space of $[X,F_Y]$ (resp. $[X,F'_Y]$) by the natural action of $[X,\Omega\hat{Y}]$ (resp. $[X,\Omega\check{Y}]$), where $F_Y$ (resp. $F'_Y$) is the homotopy fiber of the profinite completion $c_Y:Y\longrightarrow \hat{Y}$ (resp. the local expansion $e_Y:Y\longrightarrow \check{Y}$). Hence, we have the exact sequence of pointed sets
\begin{eqnarray*}
	& [X,\Omega \hat{Y}] \longrightarrow [X,F_Y] \overset{\pi}{\longrightarrow} \rmph(X,Y)\longrightarrow 0,\\
	& [X,\Omega\check{Y}] \longrightarrow [X,F'_Y] \overset{\pi'}{\longrightarrow} \rmsph(X,Y) \longrightarrow 0
\end{eqnarray*}
(see \cite[Lemma 3.5 and Corollary 5.3]{phantom}). Further, there exist (noncanonical) bijections
\begin{eqnarray*}
	& [X,F_Y] \cong \underset{i>0}{\prod} H^i (X; \pi_{i+1}(Y)\otimes \hat{\zbb}/\zbb),\\
	& [X,F'_Y] \cong \underset{i>0}{\prod} H^i(X; \pi_{i+1}(Y)\otimes \check{\zbb}/\zbb),
\end{eqnarray*}
where $\check{\zbb}$ denotes the product $\Pi_p \zbb_{(p)}$ of the $p$-localizations of $\zbb$, in which $\zbb$ is diagonally contained (see \cite[Proposition 5.4 and Remark 5.6]{phantom}).
\begin{thm}\label{thm1.7}
	Let $(f^{op},g):(X,Y)\longrightarrow (X',Y')$ be a morphism of $\ccal\wcal^{op}\times \ncal$. If $H^i(f;\pi_{i+1}(g)\otimes \mathbb{Q})$ is a monomorphism for any $i>0$,	then the following implications hold:
	\begin{eqnarray*}
		& [X', F_{Y'}] \xrightarrow[\cong]{\pi} \rmph(X',Y') \Rightarrow [X,F_Y] \xrightarrow[\cong]{\pi} \rmph(X,Y),\\
		& [X',F'_{Y'}] \xrightarrow[\cong]{\pi'} \rmsph(X',Y') \Rightarrow [X,F'_Y] \xrightarrow[\cong]{\pi'} \rmsph(X,Y).
	\end{eqnarray*}
\end{thm}
Part 2 of the following corollary is a generalization of Theorem B (see Remark \ref{rem3.3}) and Part 1 is its dual version.\vspace{2mm}
\begin{cor}\label{cor1.8}
	\begin{itemize}
		\item[{\rm (1)}] Let $f:X'\longrightarrow X$ be a map of $\ccal\wcal$, and $Y$ an object of $\ncal$. If $H_i(f;\mathbb{Q})$ is an epimorphism for $i>0$ with $\pi_{i+1}(Y)\otimes \mathbb{Q} \neq 0$, then the following implications hold:
		\begin{eqnarray*}
			& [X', F_Y] \xrightarrow[\cong]{\pi} \rmph(X',Y) \Rightarrow [X,F_Y] \xrightarrow[\cong]{\pi} \rmph(X,Y),\\
			& [X', F'_Y] \xrightarrow[\cong]{\pi'} \rmsph(X',Y) \Rightarrow [X,F'_Y] \xrightarrow[\cong]{\pi'} \rmsph(X,Y).
		\end{eqnarray*}
		\item[{\rm (2)}] Let $X$ be an object of $\ccal\wcal$ and $g:Y\longrightarrow Y'$ a map of $\ncal$. If $\pi_{i+1}(g)\otimes \mathbb{Q}$ is a monomorphism for $i>0$ with $H_i(X;\mathbb{Q})\neq 0$, then the following implications hold:
		\begin{eqnarray*}
			& [X,F_{Y'}] \xrightarrow[\cong]{\pi} \rmph(X,Y') \Rightarrow [X,F_Y] \xrightarrow[\cong]{\pi} \rmph(X,Y),\\
			& [X,F'_{Y'}] \xrightarrow[\cong]{\pi'} \rmsph(X,Y') \Rightarrow [X,F'_Y] \xrightarrow[\cong]{\pi'} \rmsph(X,Y).
		\end{eqnarray*}
	\end{itemize}
\end{cor}
Corollary \ref{cor1.8} can be applied to find various pairs $(X,Y)$ for which $\rmph(X,Y)\cong \underset{i>0}{\Pi} H^i (X;\pi_{i+1}(Y)\otimes \hat{\zbb}/\zbb)$ hold (see Corollary \ref{cor3.4}, Proposition \ref{cor1.9}, Corollary \ref{positiveintegerm} and Example \ref{newth2}).\par
The results in this section are proved in Sections \ref{section2} and 3.

\section{Proof and applications of Theorem \ref{thm1.1}} \label{section2}
In this section, we prove Theorem \ref{thm1.1} and Corollary \ref{cor1.2}, and apply Corollary \ref{cor1.2} to obtain new vanishing results for $\rmph(X,Y)$. In this and next sections, the subscript ($0$) denotes the rationalization of a nilpotent space or a nilpotent group.
\begin{proof}[\rm {Proof of Theorem \ref{thm1.1}}]
	 We can assume that the targets are in the full subcategory $1$-$\ccal\wcal$ of $\ccal\wcal$ consisting of $1$-connected $CW$-complexes (\cite[Remark 5.6]{phantom}).\par
	 Recall that there exist natural bijections
	 \begin{eqnarray*}
	 	\rmph(X,Y) &\cong& (\Omega c_{(0)})_{\sharp} [X, \Omega Y_{(0)}] \backslash [X, \Omega \hat{Y}_{(0)}]/ (\Omega\ \hat{}r)_{\sharp}[X,\Omega \hat{Y}], \\
	 	\rmsph(X,Y) &\cong& (\Omega e_{(0)})_{\sharp} [X, \Omega Y_{(0)}] \backslash [X, \Omega \check{Y}_{(0)}]/ (\Omega\ \check{ }r)_{\sharp}[X,\Omega \check{Y}],
	 \end{eqnarray*}
 	where $\hat{\,}r$ and $\check{\,}r$ denote the rationalizations of $\hat{Y}$ and $\check{Y}$ respectively (\cite[Proposition 5.7]{phantom}). Recall also that the functor $[\,\cdot\,,\Omega\,\cdot_{(0)}]:\ccal\wcal^{op} \times 1$-$\ccal\wcal \rightarrow Set$ is naturally isomorphic to the functor $\hom_\mathbb{Q}(H_\ast(\,\cdot\,;\mathbb{Q}),\pi_\ast(\Omega\,\cdot\,)\otimes \mathbb{Q})\cong \underset{i>0}{\prod}H^i(\,\cdot\,;\pi_{i+1}(\,\cdot\,)\otimes \mathbb{Q})$ (\cite[Proposition 4.1]{phantom}). Then, the result easily follows.
\end{proof}
\begin{proof}[\rm {Proof of Corollary \ref{cor1.2}}]
	The result is immediate from Theorem \ref{thm1.1}
\end{proof}
	
	In the rest of this section, we derive vanishing results for $\rmph(X,Y)$ from Corollary \ref{cor1.2}. Recall the following basic vanishing results:
	\begin{itemize}
		\item $\rmph(K,Y)=0$ for a finite complex $K$.
		\item $\rmph(\Omega K, Y)=0$ for a $1$-connected finite complex $K$ and $Y\in \ncal$.
	\end{itemize}
	The first is obvious from the definition of a phantom map. The second is a result of Iriye \cite[Theorem 1.A]{Iriye04}.\par
	The $n$-connected cover $X\langle n\rangle$ of $X$ is called $good$ if the canonical map $X\langle n\rangle \longrightarrow X$ induces a monomorphism on the rational homology.
\begin{prop}\label{prop1.4}
	Let $X$ be in $\ccal\wcal$ and $Y$ in $\ncal$, and suppose that $\rmph(X,Y)=0$. If the $n$-connected cover $X\langle n\rangle$ is good, then $\rmph(X\langle n\rangle,Y)=0$.
\begin{proof}[\rm {Proof}]
	The result follows from Corollary \ref{cor1.2}(1).
\end{proof}
\end{prop}
Recall the definition of an $H_0$-space from Section 1.\vspace{2mm}
\begin{cor}\label{cor1.5}
	\begin{itemize}
		\item[$(1)$] If $K$ is a nilpotent finite complex which is an $H_0$-space, then $\rmph(K\langle n\rangle, Y)=0$ for any $n>0$ and any $Y\in \ncal$.
		\item[$(2)$] If $K$ is a $1$-connected finite complex, then $\rmph(\Omega K\langle n \rangle, Y)=0$ for any $n>0$ and any $Y\in \ncal$.
	\end{itemize}
\begin{proof}[\rm {Proof}]
	As mentioned above, if $X$ is a finite complex or the loop space of a $1$-connected finite complex, then $\rmph(X,Y)=0$ for any $Y\in \ncal$. Note that if $X$ is an $H_0$-space, then the $n$-connected cover $X\langle n \rangle$ is good. Then, the results follow from Proposition \ref{prop1.4}.
\end{proof}
\end{cor}
\begin{rem}\label{rem1.6}{\rm
		Let $K$ be a $1$-connected finite complex with $\pi_2$ finite. Then, we can use \cite[Corollary 2.9]{phantom} along with Proposition \ref{prop1.4} to prove that $K\langle n \rangle$ is good if and only if $\rmph(K\langle n\rangle,Y)=0$ for any $Y\in \ncal$, giving an alternative proof of \cite[Theorem 6]{McGibbon98}.
}\end{rem}\vspace{2mm}
We end this section with an application of Corollary \ref{cor1.2}(2). For a $G$-space $X$, the homotopy orbit space $X_{hG}$ of $X$ (or the Borel construction on $X$) is defined by
\[
	X_{hG} = EG \times_G X,
\]
where $EG$ is the total space of the universal principal $G$-bundle (see \cite{MayG}).\vspace{2mm}
\begin{exa}\label{exa2.5}{\rm 
	Let $G$ be a compact Lie group and $X$ a finite $G$-$CW$-complex. Let $\gbb$ denote the infinite unitary group $U$ or the infinite orthogonal group O, and let $H$ be a compact Lie group which is a topological subgroup of $\gbb$. Let $k$ be a nonnegative integer. If $\pi_{i+1}(H)\otimes \qbb\longrightarrow \pi_{i+1}(\gbb)\otimes \qbb$ is a monomorphism for $i>0$ with $H_i(\Sigma^k X_{hG};\qbb)\neq 0$, then
	\[
		\rmph(X_{hG},\Omega^k \gbb/H)=0.
	\]
	\begin{proof}[\rm {Proof}]
		Recall from \cite[Section 2.3]{Gr} that
	\[
		\rmph(X_{hG},\Omega^k \gbb/H)\cong \rmph(\Sigma^k X_{hG},\gbb/H)
	\]
	and that
	\[
		\rmph(\Sigma^k X_{hG},\gbb)=0.
	\]
	Then, the result follows from Corollary \ref{cor1.2}(2).
	\end{proof}
}\end{exa}

\section{Proof and applications of Theorem \ref{thm1.7}}
In this section, we prove Theorem \ref{thm1.7} and Corollary \ref{cor1.8}, and apply Corollary \ref{cor1.8} to find various pairs $(X,Y)$ for which $\rmph(X,Y)\cong \underset{i>0}{\Pi}H^i(X;\pi_{i+1}(Y)\otimes \hat{\zbb}/\zbb)$ hold.\par
For the proof of Theorem \ref{thm1.7}, we need the following lemma.
\begin{lem}\label{lem3.1}
	Let $G$ and $H$ be groups, and $G'$ and $H'$ subgroups of $G$ and $H$ respectively. Consider a morphism of short exact sequences of pointed sets
	\begin{center}	
		\begin{tikzcd}
		\ast \arrow[r] & G' \arrow[r] \arrow[d, hook'] & G \arrow[r] \arrow[d, hook'] & G/G' \arrow[r] \arrow[d] & \ast \\
		\ast \arrow[r] & H' \arrow[r] & H \arrow[r] & H/H' \arrow[r] & \ast
		\end{tikzcd}
	\end{center}
	where left and middle vertical arrows are monomorphisms of groups. Then, the right vertical arrow is injective if and only if the left square is a pullback diagram in Set.
	\begin{proof}[\rm {Proof}]
		$(\Leftarrow)$ Let $\bar{g_1}$ and $\bar{g_2}$ be elements of $G/G'$ such that $\bar{g_1}=\bar{g_2}$ in $H/H'$. Then, $g_1h'=g_2$ for some $h'\in H'$. Since $h'=g_1^{-1}g_2 \in G$, $h' \in H' \cap G = G'$, and hence $\bar{g_1}=\bar{g_2}$ in $G/G'$.\\
		$(\Rightarrow)$ Suppose that the left square is not a pullback diagram in $Set$. Then, we can choose an element $h'\in H'\cap G\backslash G'$. Since $\bar{h'}\neq \bar{e}$ in $G/G'$ and $\bar{h'}=\bar{e}$ in $H/H'$, the right vertical arrow is not injective.
	\end{proof}
\end{lem}

\begin{proof}[\rm {Proof of Theorem \ref{thm1.7}}]
	We prove the first implication; the proof of the second implication is similar.\par
	By \cite[Remark 5.6]{phantom}, we assume that $Y$ and $Y'$ are $1$-connected.\par
	For a $1$-connected $CW$-complex $Z$ of finite type, we have the fibration sequence
	\[
		F_Z \longrightarrow Z \overset{c_Z}{\longrightarrow} \hat{Z}.
	\]
	Since $F_Z$ is a simple rational space (\cite[Proposition 5.4]{phantom}), we obtain from this the fibration sequence
	\[
		\Omega Z_{(0)} \longrightarrow \Omega \hat{Z}_{(0)} \longrightarrow F_Z.
	\]
	Since this is a trivial principal fibration (\cite[Proposition 3.1]{phantom}), the map $[A,\Omega Z_{(0)}]\longrightarrow [A,\Omega \hat{Z}_{(0)}]$ is a monomorphism of groups and the map $[A,\Omega \hat{Z}_{(0)}]\longrightarrow [A,F_Y]$ can be identified with the quotient map for the action of $[A,\Omega Z_{(0)}]$ for any $A\in \ccal\wcal$.\par
	Thus, we consider the morphism of short exact sequences of pointed sets
	\[
	\begin{tikzcd}
			\ast \arrow{r} & \,[X,\Omega Y_{(0)}] \arrow[r] \arrow[swap]{d}{(f^{op},\Omega g_{(0)})_\sharp} & \,[X, \Omega\hat{Y}_{(0)}] \arrow[r] \arrow[swap]{d}{(f^{op},\Omega\hat{g}_{(0)})_\sharp}& \,[X, F_Y] \arrow[r] \arrow[swap]{d}{(f^{op},F_g)_\sharp} & \ast\\
			\ast \arrow[r] & \,[X',\Omega Y'_{(0)}] \arrow[r] & \,[X',\Omega\hat{Y'}_{(0)}] \arrow[r] & \,[X',F_{Y'}] \arrow[r] & \ast \tag{3.1}
	\end{tikzcd}
	\]
	Recall that the functor $[\,\cdot\,,\Omega\,\cdot_{(0)}]:\ccal\wcal^{op} \times 1$-$\ccal\wcal \rightarrow Set$ is naturally isomorphic to the functor $\underset{i>0}{\prod}H^i(\,\cdot\,;\pi_{i+1}(\,\cdot\,)\otimes \mathbb{Q})$ (\cite[Proposition 4.1]{phantom}) and identify the left square in (3.1) with the square
	\[
	\begin{tikzcd}
		\underset{i>0}{\prod}H^i(X; \pi_{i+1}(Y)_{(0)}) \arrow{r} \arrow[swap]{d}{\underset{i>0}{\prod}H^i(f;\pi_{i+1}(g)_{(0)})} & \underset{i>0}{\prod}H^i (X;\pi_{i+1}(Y)\otimes \hat{\zbb}_{(0)}) \arrow{d}{\underset{i>0}{\prod}H^i(f;\pi_{i+1}(g)\otimes \hat{\zbb}_{(0)})}\\
		\underset{i>0}{\prod}H^i(X';\pi_{i+1}(Y')_{(0)}) \arrow{r} & \underset{i>0}{\prod}H^i(X';\pi_{i+1}(Y')\otimes \hat{\zbb}_{(0)}) \tag{3.2}
	\end{tikzcd}
	\]
	in the category of sets. Then, we see from the assumption that the two vertical arrows in (3.2), and hence the left and middle vertical arrows in (3.1) are injective.\par
	Now, consider the morphism of short exact sequences of $\qbb$-modules
	\begin{center}
	\begin{tikzcd}
		0 \arrow{r} & \underset{i>0}{\prod}H^i(X; \pi_{i+1}(Y)_{(0)}) \arrow{r} \arrow[swap]{d}{\underset{i>0}{\prod}H^i(f;\pi_{i+1}(g)_{(0)})} & \underset{i>0}{\prod}H^i (X;\pi_{i+1}(Y)\otimes \hat{\zbb}_{(0)})\arrow{d}{\underset{i>0}{\prod}H^i(f;\pi_{i+1}(g)\otimes \hat{\zbb}_{(0)})} \arrow{r} & \underset{i>0}{\prod}H^i(X;\pi_{i+1}(Y)\otimes \hat{\zbb}/\zbb) \arrow{r} \arrow{d}{\underset{i>0}{\prod}H^i(f;\pi_{i+1}(g)\otimes \hat{\zbb}/\zbb)} & 0\\
		0 \arrow{r} & \underset{i>0}{\prod}H^i(X';\pi_{i+1}(Y')_{(0)}) \arrow{r} & \underset{i>0}{\prod}H^i(X';\pi_{i+1}(Y')\otimes \hat{\zbb}_{(0)}) \arrow{r} & \underset{i>0}{\prod}H^i(X';\pi_{i+1}(Y')\otimes \hat{\zbb}/\zbb) \arrow{r} & 0
	\end{tikzcd}
	\end{center}
	which extends (3.2). Since the right vertical arrow is injective by the assumption, the left square is a pullback diagram in $Set$ (Lemma \ref{lem3.1}). Thus, the right vertical arrow in (3.1) is injective by Lemma \ref{lem3.1}.\par
	Next consider the commutative diagram
	\[
		\begin{tikzcd}
		\ [X,F_Y] \arrow{r}{(f^{op},F_g)_\sharp} \arrow[d, "\pi"'] & \ [X',F_{Y'}]\arrow[d, "\pi"] \\
		\rmph(X,Y) \arrow{r}{(f^{op},g)_\sharp} & \rmph(X',Y')
		\end{tikzcd}
	\]
	and recall that the vertical arrows are quotient maps. Since the upper horizontal arrow $(f^{op},F_g)_\sharp$ is injective, the implication in question is obvious.
\end{proof}
\begin{proof}[\rm {Proof of Corollary \ref{cor1.8}}]
	The result is immediate from Theorem \ref{thm1.7}.	
\end{proof}
\begin{rem}\label{rem4.2}{\rm
	For $(X,Y)\in \ccal\wcal^{\rm op} \times \ncal$, there exist no natural bijections
	\begin{eqnarray*}
		& [X,F_Y] \cong \underset{i>0}{\prod} H^i (X; \pi_{i+1}(Y)\otimes \hat{\zbb}/\zbb),\\
		& [X,F'_Y] \cong \underset{i>0}{\prod} H^i(X; \pi_{i+1}(Y)\otimes \check{\zbb}/\zbb).
	\end{eqnarray*}
	However, if we restrict ourselves to the subclass $\qcal$ of $\ccal\wcal^{\rm op} \times \ncal$, then the above bijections can be taken to be natural ones (\cite[Proposition 5.10]{phantom}), which induce natural divisible abelian group structures on $\rmph(X,Y)$ and $\rmsph(X,Y)$ (\cite[Theorem 2.3]{phantom}). Here, the subclass $\qcal$ is defined by the following condition:
	\begin{itemize}
		\rmitem{Q} For each pair $i,j>0$, the rational cup product
		\[\ \cup:\ H^{i}(X;{\mathbb Q}) \otimes H^{j}(X;{\mathbb Q})\longrightarrow H^{i+j}(X;{\mathbb Q})\]
		or the rational Whitehead product
		\[[\ ,\ ]:\ (\pi_{i+1}(Y)\otimes {\mathbb Q}) \otimes (\pi_{j+1}(Y)\otimes {\mathbb Q}) \longrightarrow \pi_{i+j+1}(Y)\otimes {\mathbb Q} \]
		is trivial. 
	\end{itemize}
}\end{rem}\vspace{2mm}

		Let us recall the generalizations of Miller's theorem \cite{Miller84} and Anderson-Hodgkin's theorem \cite{Anderson68} to obtain many pairs $(X, Y')$ with $[X,F_{Y'}] \xrightarrow[\cong]{\pi} \rmph(X,Y')$. A space whose $i^{\rm th}$ homotopy group is zero for $i \leq n$ and locally finite for $i = n + 1$ is said to be {\it $n\frac{1}{2}$-connected}.
Define the classes $\mathcal A$, $\mathcal B$, $\mathcal A'$ and
$\mathcal B'$ by
\begin{enumerate}
	\item[$\mathcal A$ =] the class of $\frac{1}{2}$-connected Postnikov spaces, the classifying spaces of compact Lie groups, 
	$\frac{1}{2}$-connected infinite loop spaces and their iterated suspensions.
	
	\item[$\mathcal B$ =] the class of nilpotent finite complexes, the classifying spaces of compact Lie groups and their iterated loop spaces.
	
	\item[$\mathcal A'$ =] the class of ${1}\frac{1}{2}$-connected
	Postnikov spaces of finite type and their iterated suspensions.
	
	\item[$\mathcal B'$ =] the class of $BU$, $BO$, $BSp$, $BSO$, $U/Sp$,
	$Sp/U$, $SO/U$, $U/SO$, and their iterated loop spaces.
\end{enumerate}
If $(X, Y')$ is in $\acal \times \bcal$ or $\acal^{'} \times \bcal^{'}$, then ${[X, \Omega \hat{Y'}] = 0}$ (\cite[Corollary 6.4]{phantom}), and hence $[X,F_{Y'}] \xrightarrow[\cong]{\pi} \rmph(X,Y')$ and $[X,F'_{Y'}] \xrightarrow[\cong]{\pi'} \rmsph(X,Y')$ (\cite[Propositions 6.1 and 2.5]{phantom}).\vspace{2mm}
\begin{rem}\label{rem3.3}{\rm
	Note that $K(\zbb,n)$ is in $\acal'$ for $n\geq 3$ and that $U$ and $BU$ are in $\bcal'$. Also note that a space $Y$ as in Theorem B has trivial rational Whitehead products. Then, we see that Corollary \ref{cor1.8}(2) implies Theorem B (see Remark \ref{rem4.2}).
}\end{rem}\vspace{2mm}
Part 1 of the following corollary is a direct generalization of \cite[Corollary 2 in Section 2]{Meier78}. For a connected $CW$-complex $K$, $Q(K)$ denotes the infinite loop space defined by $Q(K)=\underset{n}{\rm colim}\,\Omega^n \Sigma^n K$.
	\begin{cor}\label{cor3.4}
		\begin{itemize}
			\item[{\rm (1)}] Let $X$ be a space in $\acal'$. Then, there exist natural isomorphisms of abelian groups
			\begin{eqnarray*}
				\rmph(X,QS^p)\cong H^{p-1}(X;\hat{\zbb}/\zbb),\\
				\rmsph(X,QS^p)\cong H^{p-1}(X;\check{\zbb}/\zbb).
			\end{eqnarray*}
			\item[{\rm (2)}] Let $K$ be a $2\frac{1}{2}$-connected finite complex. If $n\geq \dim K$, then there exist natural isomorphisms of abelian groups
			\begin{eqnarray*}
				\rmph(K\langle n\rangle, QS^p)\cong H^{p-1}(K\langle n\rangle; \hat{\zbb}/\zbb),\\
				\rmsph(K\langle n\rangle, QS^p)\cong H^{p-1}(K\langle n\rangle; \check{\zbb}/\zbb).
			\end{eqnarray*}
		\end{itemize}
	\begin{proof}[\rm {Proof}]
		See Remark \ref{rem4.2} for the natural abelian group structures on $\rmph(A,QS^p)$ and $\rmsph(A,QS^p)$.
		\begin{itemize}
			\item[{\rm (1)}] $QS^p$ admits a map to $U$ (resp. $BU$) which induces monomorphisms on the rational homotopy groups for $p$ odd (resp. even); see \cite[p. 476]{Meier78}. Thus, the result follows from \cite[Proposition 2.5]{phantom} and Corollary \ref{cor1.8}(2).
			\item[{\rm (2)}] We see from \cite[Corollary 2.9]{phantom} that
			\[
				[K\langle n\rangle, F_{Y'}] \overset{\pi}{\longrightarrow} \rmph(K\langle n \rangle,Y') \text{ and } [K\langle n\rangle, F'_{Y'}] \overset{\pi'}{\longrightarrow} \rmsph(K\langle n\rangle,Y')
			\]
			are bijective for $Y'\in \bcal'$. Thus, the result follows from Corollary \ref{cor1.8}(2) (see the proof of Part 1).\qedhere
		\end{itemize}
	\end{proof}
	\end{cor}
\begin{prop}\label{cor1.9}
	Let $X$ be in $\ccal\wcal$ and $Y'$ in $\ncal$. Then, the following implications hold for any $m>0$;
	\begin{eqnarray*}
		& [X,F_{Y'}] \xrightarrow[\cong]{\pi} \rmph(X,Y') \Rightarrow [X,F_{Y'\langle m \rangle}] \xrightarrow[\cong]{\pi} \rmph(X,Y'\langle m \rangle),\\
		& [X,F'_{Y'}] \xrightarrow[\cong]{\pi} \rmsph(X,Y') \Rightarrow [X,F'_{Y'\langle m \rangle}] \xrightarrow[\cong]{\pi} \rmsph(X,Y'\langle m \rangle).
	\end{eqnarray*}
\begin{proof}[\rm {Proof}]
	The result follows from Corollary \ref{cor1.8}(2).
\end{proof}
\end{prop}
\begin{cor}\label{positiveintegerm}
	Let $(X,Y')$ be in $\mathcal{A}\times\mathcal{B}$ or $\mathcal{A}'\times\mathcal{B}'$, and let $m$ be a positive integer. Then there exist (noncanonical) bijections
	\begin{eqnarray*}
		\rmph(X, Y'\langle m\rangle) &\cong& \underset{i>0}{\prod}\ H^i (X; \pi_{i+1}(Y'\langle m\rangle) \otimes \hat{\mathbb Z} / \mathbb Z), \\
		\rmsph(X, Y'\langle m\rangle ) &\cong& \underset{i>0}{\prod}\ H^i (X; \pi_{i+1}(Y'\langle m\rangle) \otimes \check{\mathbb Z}/{\mathbb Z}).
	\end{eqnarray*}
\begin{proof}[\rm {Proof}]
	The result follows from Proposition \ref{cor1.9}.
\end{proof}
\end{cor}
	Corollary \ref{positiveintegerm}, and hence Proposition \ref{cor1.9} can be regarded as a generalization of \cite[Corollary 1.2]{Gr} (see Remark \ref{rem4.2}).\par
	We end this section with an application of Corollary \ref{cor1.8}(1). The Grassmannians $G_n(\mathbb{F})$ and $G_{\infty}(\mathbb{F})$ are defined by $G_n(\mathbb{F})=\underset{\longrightarrow}{\lim}_m\ G_{n,m}(\mathbb{F})$ and $G_{\infty}(\mathbb{F})=\underset{\longrightarrow}{\lim}_n\ G_n(\mathbb{F})$ for $\mathbb{F}=\mathbb{R},\mathbb{C},\mathbb{H}$, where the finite Grassmannian $G_{n,m}(\mathbb{F})$ is the space of $n$-dimensional subspaces in $\mathbb{F}^{n+m}$.\vspace{2mm}

\begin{exa}\label{newth2}
	{\rm Let $G_{n}(\fbb)/G_{n', m'}(\fbb)$ be the quotient complex of $G_{n}(\fbb)$ by $G_{n', m'}(\fbb)$ $(n', m' < \infty,\ n' \leq n \leq \infty)$. Let $Y$ be a space in $\bcal$. If $\mathbb{F}=\mathbb{C}$ or $\mathbb{H}$, then there exist (noncanonical) bijections}
	\begin{eqnarray*}
		\rmph(G_n(\mathbb{F})/G_{n',m'}(\mathbb{F}), Y) & \cong & \underset{i>0}{\prod}\, H^i(G_n(\mathbb{F})/G_{n',m'}(\mathbb{F}) ; \pi_{i+1}(Y) \otimes \hat{\mathbb Z} / \mathbb Z), \\
		\rmsph(G_n(\mathbb{F})/G_{n',m'}(\mathbb{F}), Y) & \cong & \underset{i>0}{\prod}\, H^i(G_n(\mathbb{F})/G_{n',m'}(\mathbb{F}) ; \pi_{i+1}(Y) \otimes \check{\mathbb Z} / \mathbb Z).
	\end{eqnarray*}
\end{exa}
\begin{proof}[\rm {Proof}]
	Note that the quotient map $G_n(\mathbb{F})\longrightarrow G_n(\mathbb{F})/G_{n',m'}(\mathbb{F})$ induces an epimorphism on the rational homology. Since $G_n(\mathbb{F})$ is in $\acal$, the maps
	\[
		[G_n(\mathbb{F}),F_Y] \overset{\pi}{\longrightarrow} \rmph(G_n(\mathbb{F}),Y) \text{  and  } [G_n(\mathbb{F}),F'_Y] \overset{\pi'}{\longrightarrow} \rmsph(G_n(\mathbb{F}),Y)
	\]
	are bijective (\cite[Proposition 2.5]{phantom}). Thus, the result follows from Corollary \ref{cor1.8}(1).
\end{proof}
	If $(G_n(\mathbb{F})/G_{n',m'}(\mathbb{F}), Y)$ is in $\qcal$, then we can obtain a more precise result (see \cite[Example 6.6]{phantom}).

\end{document}